\documentstyle[pb-diagram]{article}
\newcommand{\kitem}{\begin{itemize}\vspace{-2ex}}
\newcommand{\kenditem}{\vspace{-1ex}\end{itemize}}

\newcommand{\bean}{\[\begin{array}{rcl}}
\newcommand{\eean}{\end{array}\]}
\newcommand{\N}{{I\!\!N}}
\newcommand{\R}{{I\!\!R}}
\newcommand{\Z}{{Z\!\!\!Z}}

\newcommand{\veee}{{\scriptscriptstyle\vee}}
\newcommand{\Spec}{\mbox{\rm Spec}}

\newcommand{\tot}{\mbox{\rm tot}}
\newcommand{\innt}{\mbox{\rm int}}
\newcommand{\kar}{\mbox{\rm char}}
\newcommand{\orb}{\mbox{\rm orb}}
\newcommand{\spann}{\mbox{\rm span}}

\newcommand{\ko}{\overline}
\newcommand{\ku}{\underline}
\newcommand{\im}{\mbox{\rm im}\,}

\newcommand{\kss}{\scriptscriptstyle}
\newcommand{\kd}{\displaystyle}
\newcommand{\kf}{\footnotesize}

\newcommand{\gExt}{\mbox{\rm Ext}}

\newcommand{\gHom}{\mbox{\rm Hom}}

\newcommand{\an}{m}
\newcommand{\ai}{i}
\newcommand{\aj}{j}
\newcommand{\ak}{k}
\newcommand{\al}{l}

\newcommand{\HA}{H\!A}
\newcommand{\kL}{\Lambda}
\newcommand{\bHom}{\ko{\mbox{\rm Hom}}\,}
\newcommand{\bL}{\ko{L}}
\newcommand{\kq}{q}

\newcounter{Abschnitt}[section]
\newcommand{\neu}[1]{\protect\refstepcounter{Abschnitt}\protect
   \vspace{1ex}\label{#1}
   {\bf (\protect\arabic{section}.\protect\arabic{Abschnitt})}
                     $\qquad$}
\newcommand{\zitat}[2]{(\protect\ref{#1}.\protect\ref{#1-#2})}

\begin{document}
\title{Andr\'{e}-Quillen cohomology of monoid algebras}

\author{Klaus~Altmann\qquad Arne~B.~Sletsj\o{}e}
\date{}
\maketitle


\begin{abstract}
We compute the Andr\'{e}-Quillen (or Harrison)
cohomology of an affine toric variety. The best results are
obtained either in the general case for the first three cohomology groups, or in the
case of
isolated singularities for all cohomology groups, respectively.
\end{abstract}

%
%
\section{Introduction}\label{Int}


\neu{Int-1}
Let $k$ be a field. For any finitely generated $k$-algebra $A$ 
the so-called cotangent complex yielding the
Andr\'{e}-Quillen cohomology $T^n_A=T^n(A,A;k)$ ($n\geq 0$) may be defined. 
The first three
of these $A$-modules are important for the deformation theory of $A$ or its geometric
equivalent $\Spec A$: $T^1_A$ equals the set of infinitesimal deformations, $T^0_A$
describes their automorphisms, and $T^2_A$ contains the obstructions for lifting
infinitesimal deformations to larger base spaces. 
Apart from occurring in long exact sequences 
no meaning of the higher cohomology groups seems to be known
when studying the deformation theory of closed subsets of $\Spec A$.
A very readable reference for the definition of Andr\'{e}-Quillen cohomology 
and its
relations to Hochschild and Harrison cohomology is Loday's book \cite{Loday}. For
applications in deformation theory see for instance \cite{LaudalSLN}, \cite{Pa},
or the summary of the
properties one has to know (without proofs) in the first section of \cite{BeCh}.
\par


\neu{Int-2}
For smooth $k$-algebras $A$, the higher $T^n_A$ (i.e.\ $n\geq 1$) vanish. For
complete intersections the situation is still easy; $T^0_A$ and $T^1_A$ are well
understood, and the remaining cohomology groups vanish. As far as we know, only 
few further
examples exist where the cotangent complex or at least the Andr\'{e}-Quillen
cohomology groups are known. Palamodov has told us that he has computed (unpublished) 
the cotangent complex of an embedded point on a line; it turned out that the
Poincar\'{e} series $T(s):=\sum_{n\geq 0}(\dim_k T^n)\cdot s^n$ of this singularity
is a rational function. It would be interesting to know whether this is always the
case for isolated singularities.\\
The result of the present paper is a spectral sequence converging to the Harrison (or
Andr\'{e}-Quillen) cohomology for affine toric varieties. 
In the case of an isolated singularity, this spectral
sequence degenerates; this leads to a down to earth description of the modules $T_A^n$.
Moreover, in the general case, the information is still sufficient to determine
$T_A^0$, $T_A^1$, and $T^2_A$. That is, we use the methods of \cite{Sl} to obtain
straight formulas generalizing part of the results of \cite{T2}.
\par


\neu{Int-3}
The paper is organized as follows:
We begin in \S \ref{iHc} with fixing notation and recalling those facts of 
\cite{Sl} that will be used in the following. The sections \S \ref{sps} and
\S \ref{E1l} contain the main theorems of the paper: First, we state our spectral
sequence; then, in
\S \ref{E1l}, we calculate some of the $E_1$-terms and show the vanishing of others.
Finally, in \S \ref{app}, we present the resulting $T^n_A$-formulas promised
before.
\par


\neu{Int-4} {\em Acknowledgements:}
The first author became interested in describing the algebra cohomology of 
monoid algebras during his
visit at the Senter for h\o{}yere studier (Academy of Sciences) in Oslo. 
He liked the
warm and stimulating atmosphere very much, and would like to thank those who made
it possible for him to stay there. Moreover, we are grateful to R.-O.~Buchweitz,
J.A.~Christophersen, and A.~Laudal for many helpful discussions.
\par

%
%
\section{Inhomogeneous Harrison cohomology}\label{iHc}


\neu{iHc-1}
Let $M$, $N$ be mutually dual, finitely generated, free Abelian groups; 
we denote by $M_{\R}$, $N_{\R}$ the associated real vector spaces obtained via 
base change with $\R$.
Assume we are given a rational, polyhedral cone 
$\sigma=\langle a^1,\dots,a^\an\rangle \subseteq N_{\R}$
with apex in $0$ and with $a^1,\dots,a^\an\in N$ denoting its {\em primitive} fundamental
generators (i.e.\ none of the $a^\ai$ is a proper multiple of an element of $N$).
We define the dual cone
$\sigma^{\veee}:= \{ r\in M_{\R}\,|\; \langle \sigma,\,r\rangle \geq 0\} 
\subseteq M_{\R}$ and denote by $\kL:=\sigma^{\veee}\cap M$ the resulting monoid
of lattice points.\\
The corresponding monoid algebra $A:= k[\kL]$ will be the object of the upcoming
investigations. It is the ring of regular
functions on the toric variety $Y_\sigma= \Spec A$ associated to $\sigma$.
The ring $A$ itself as well as most of its important modules (such as $T^n_A$)
admit an $M$-(multi)grading. It is this grading which will make computations possible.
For general facts concerning toric varieties see for instance \cite{Oda}.
\par


\neu{iHc-2}
The following definitions are taken from \cite{Sl}, \S 2. 
We note that the original requirement that the monoids involved have no
non-trivial subgroups is unnecessary for our purposes.
\par

{\bf Definition:} 
{\em
$L\subseteq \kL$ is said to be monoid-like if for all elements 
$\lambda_1,\lambda_2\in L$ the relation $\lambda_1-\lambda_2\in 
\kL_+:=\kL\setminus\{0\}$ implies $\lambda_1-\lambda_2\in L$.\\
Moreover, a subset $L_0\subseteq L$ of a monoid-like set is called full if
$(L_0 + \kL)\cap L =L_0$.
}
\par

For any subset $P\subseteq \kL$ and $n\geq 1$ we introduce 
$S_n(P):= \{(\lambda_1,\dots,\lambda_n)\in P^n\,|\; \sum_v \lambda_v \in P\}$.
If $L_0\subset L$ are as in the previous definition, then this gives rise to the
following set:
\[
C^n(L\setminus L_0,L;k):= \{\varphi: S_n(L)\to k\,|\;
\varphi \mbox{ is shuffle invariant and vanishes on } S_n(L\setminus L_0)\}\,.
\]
($\varphi$ is said to be shuffle invariant if it is invariant under 
$\sum_\pi \mbox{sgn}(\pi)\cdot \pi$ where $\pi$ runs through all shuffles of the set
$\{1,\dots,n\}$.) These $k$-vector spaces turn into a complex with the
differential
$\delta^n: C^{n-1}(L\setminus L_0,L;k)\to C^n(L\setminus L_0,L;k)$ defined via
\[
(\delta^n\varphi)(\lambda_1,\dots,\lambda_n):= \varphi(\lambda_2,\dots,\lambda_n)
+ \sum_{v=1}^{n-1} (-1)^v \varphi(\lambda_1,\dots,\lambda_v+\lambda_{v+1},\dots,
\lambda_n) + (-1)^n \varphi(\lambda_1,\dots,\lambda_{n-1})\,.
\]
{\bf Definition:}
{\em
The $k$-vector space
$\HA^n(L\setminus L_0,L;k):=H^n\big(C^{\kss\bullet}(L\setminus L_0,L;k)\big)$
is called the
$n$-th inhomogeneous Harrison cohomology of the pair $(L,L_0)$.
}
\par


\neu{iHc-3}
{\bf Theorem} (\cite{Sl}):
{\em
Let $R\in M$. Then, defining $\kL_+:=\kL\setminus\{0\}$,
the homogeneous part of $T^n_A$ in degree $-R$ equals
\[
T^n_A(-R)=
\HA^{n+1}\big(\kL_+\setminus (R+\kL),\,\kL_+;\,k\big)
\quad\mbox{ for } n\geq 0.
\vspace{-3ex}
\]
}
\par

The proof of the theorem is spread throughout the first two sections of \cite{Sl}:
First, in (1.13), (1.14), the calculation of $T^n_A$ has been reduced down to
the monoid level.
Then, Proposition (2.9) shows that the homogeneous pieces $T^n_A(-R)$ equal the
so-called graded Harrison cohomology groups $\mbox{Harr}^{n+1,-R}(\kL,k[\kL])$,
and, finally, Theorem (2.10) states the above result.
\par

%
%
\section{The spectral sequence}\label{sps}


\neu{sps-1}
With the notation of \zitat{iHc}{1} we define for any face $\tau\leq\sigma$ and any
degree $R\in M$ the monoid-like set
\[
K_\tau^R:= \kL_+ \cap \big(R-\innt\, \tau^\veee\big)\,.
\]
These sets admit the following elementary properties:
\kitem
\item[(i)]
$K_0^R=\kL_+$, and $K_\ai^R:=K_{a^\ai}^R=\{r\in\kL_+\,|\; \langle a^\ai,r\rangle
<\langle a^\ai,R\rangle\}$ with $\ai=1,\dots,\an$.
\item[(ii)]
For $\tau\neq 0$ the equality $K_\tau^R=\bigcap_{a^\ai\in\tau}K_\ai^R$ holds. 
Moreover, if
$\sigma$ is a top-dimensional cone, $K_\sigma^R=\kL_+ \cap \big(R-\innt\,
\sigma^\veee\big)$ is a (diamond shaped) finite set.
\item[(iii)]
$\kL_+\setminus (R+\kL)=\bigcup_{\ai=1}^\an K_\ai^R$.
\kenditem


\neu{sps-2}
Let us fix an element $R\in M$. The dependence of the sets $K^R_\tau$ on
$\tau$ is a contravariant functor. This gives rise to the complexes
$C^q(K^R_{\kss \bullet};k)$ ($q\geq 1$) defined as
\[
C^q(K^R_p;k):= \oplus_{\,\tau\leq\sigma,\, \mbox{\kf dim}\hspace{0.1em}
               \tau=p\,} C^q(K^R_\tau;k)\qquad (0\leq p\leq \dim \sigma)
\]
with $C^q(K^R_\tau;k):= C^q(K^R_\tau,K^R_\tau;k)$ and
the obvious differentials $d^p:C^q(K^R_{p-1};k)\to C^q(K^R_p;k)$.
(One has to use the maps $C^q(K^R_\tau;k)\to C^q(K^R_{\tau^\prime};k)$
for any pair $\tau\leq\tau^\prime$ of $(p-1)$- and $p$-dimensional faces,
respectively. The only problem might be the sign; it arises from comparison
of the (pre-fixed) orientations of $\tau$ and $\tau^\prime$.)
Our complex begins as
\[
0\to C^q(\kL_+;k)\to \oplus_{\ai=1}^\an C^q(K_\ai^R;k)\to
\oplus_{\langle a^\ai,a^\aj\rangle\leq\sigma} C^q(K^R_{\ai\aj};k)\to
\oplus_{\,\mbox{\kf dim}\hspace{0.1em} \tau=3\,} C^q(K^R_\tau;k)\to
\cdots\,.
\]
\par

{\bf Lemma:}
{\em
The canonical $k$-linear map
$C^q(\kL_+\setminus (R+\kL),\,\kL_+;\,k)\to C^q(K^R_{\kss \bullet};k)$ is a
quasiisomorphism, i.e.\ a resolution of the first vector space.
}
\par

{\bf Proof:}
For an $r\in\kL_+\subseteq M$ we define the $k$-vector space
\[
V^q(r):= \big\{\varphi:\{\ku{\lambda}\in\kL_+^q\,|\,
\mbox{$\sum_v$}\lambda_v=r\}
\to k\,\big|\; \varphi \mbox{ is shuffle invariant}\big\}\,.
\]
Then, our complex splits into a direct product over $r\in\kL_+$. Its
homogeneous factors equal
\[
0\to
V^q(r) \to
V^q(r)^{\{\ai\,|\; r\in K^R_\ai\}}\to
V^q(r)^{\{\tau\leq\sigma\,|\;
\mbox{\kf dim}\hspace{0.1em}\tau=2;\,r\in K^R_\tau\}}\to
V^q(r)^{\{\tau\leq\sigma\,|\;
\mbox{\kf dim}\hspace{0.1em}\tau=3;\,r\in K^R_\tau\}}\to\cdots\,.
\]
On the other hand, denoting by $H^+_{r,R}$ the halfspace
$H^+_{r,R} := \{ a\in N_{\R}\,|\; \langle a,r\rangle < \langle a, R\rangle\}
\subseteq N_{\R}$, the relation $r\in K^R_\tau$ is equivalent to
$\tau \setminus \{0\} \subseteq H^+_{r,R}$. Hence, the complex for
computing the reduced cohomology of the topological space
$\bigcup_{\tau\setminus\{0\}\subseteq H^+_{r,R}}
\!\big(\tau\setminus \{0\}\big)\subseteq \sigma$ equals
\[
0\to k \to
k^{\{\ai\,|\; r\in K^R_\ai\}}\to
k^{\{\tau\leq\sigma\,|\;
\mbox{\kf dim}\hspace{0.1em}\tau=2;\,r\in K^R_\tau\}}\to
k^{\{\tau\leq\sigma\,|\;
\mbox{\kf dim}\hspace{0.1em}\tau=3;\,r\in K^R_\tau\}}\to\cdots
\]
if $\sigma\cap H^+_{r,R}\neq\emptyset$, i.e.\ $r\in\bigcup_\ai K^R_\ai$;
it is trivial otherwise.
Since $\bigcup_{\tau\setminus\{0\}\subseteq H^+_{r,R}}
\!\big(\tau\setminus \{0\}\big)$ is contractible, this complex is always
exact.
Thus, $C^q(K^R_{\kss \bullet};k)=
\prod_{r\in\kL_+} V^q(r)^{\{\tau\leq\sigma\,|\;
\mbox{\kf dim}\hspace{0.1em}\tau={\kss \bullet};\,r\in K^R_\tau\}}$
has
$\prod_{r\in\kL_+\setminus(\cup_\ai K^R_\ai)}V^q(r)=
C^q(\kL_+\setminus (R+\kL),\,\kL_+;\,k)$
as cohomology in $0$, and it is exact elsewhere.
\hfill$\Box$
\par


\neu{sps-3}
Combining the differentials $d^p$ from \zitat{sps}{2} and
$\delta^q$ from \zitat{iHc}{2}, we obtain a double complex
$C^{\kss \bullet}(K^R_{\kss \bullet};k)$ ($0\leq p\leq \dim\sigma;\;q\geq 1$).
\par

{\bf Theorem:}
{\em
The Andr\'{e}-Quillen cohomology of $A=k[\kL]$ equals the cohomology of the
total complex, that is
\[
T^n_A(-R)= H^{n+1}\big(\tot^{\kss \bullet}
\big[C^{\kss \bullet}(K^R_{\kss \bullet};k)\big]\big)
\quad\mbox{for } n\geq 0.
\]
Moreover, given an element $s\in\kL$, the multiplication
$[\cdot x^s]:T^n_A(-R)\to T^n_A(-R+s)$ is given by the
homomorphism $C^{\kss \bullet}(K^R_{\kss \bullet};k)\to
C^{\kss \bullet}(K^{R-s}_{\kss \bullet};k)$
assigned to the inclusions $K^{R-s}_\tau\subseteq K^R_\tau$.
}
\par

{\bf Proof:}
The first part is a straightforward consequence of Theorem \zitat{iHc}{3} and
the previous lemma -- just use the corresponding spectral sequence
of the double complex. For the $A$-module structure of $T^n_A$,
one easily observes
from the proofs in \cite{Sl} that the multiplication with $x^s$ arises
from the complex homomorphism
$C^{\kss\bullet}(\kL_+\setminus (R+\kL),\,\kL_+;\,k) \to
C^{\kss\bullet}(\kL_+\setminus ([R-s]+\kL),\,\kL_+;\,k)$
provided by the inclusion $R+\kL\subseteq [R-s]+\kL$.
\hfill$\Box$
\par


\neu{sps-4}
\hspace*{-0.5em}{\bf Corollary:}
{\em
There is a spectral sequence
$E_1^{p,q}=\oplus_{\mbox{\kf dim}\hspace{0.1em}\tau=p} \HA^q(K^R_\tau;k)
\Longrightarrow T_A^{p+q-1}(-R)$.
}
\par

{\bf Proof:}
This is the other spectral sequence associated to the double complex
$C^{\kss \bullet}(K^R_{\kss \bullet};k)$.
\hfill$\Box$
\par

%
%
\section{The $E_1$-level}\label{E1l}


\neu{E1l-1}
{\bf Definition:}
Let $K\subseteq M$ be an arbitrary subset of the lattice $M$.
A function $f:K\to k$
is called {\em quasilinear} if $f(r)+f(s)=f(r+s)$ for any $r$ and $s$
with $r, s, r+s\in K$. The vector space of quasilinear functions is denoted by
$\bHom(K,k)$.
\par

Recalling the differential $\delta^2:C^1(K^R_\tau)\to C^2(K^R_\tau)$
from \zitat{iHc}{2} shows that the $E_1^{{\kss \bullet},1}$-summands
$\HA^1(K^R_\tau;k)$ equal $\bHom(K^R_\tau,k)$.
\par


\neu{E1l-2}
The orbits of the torus acting on $Y_\sigma$ are parametrized by the faces
of $\sigma$; the singular locus of $Y_\sigma$ is the disjoint union of some of these
orbits. We call a face $\tau\leq\sigma$ smooth if our toric variety is smooth along
$\orb(\tau)$. It is one of the basic facts that smooth faces are characterized by
being generated from a part of a $\Z$-basis of $N$. In particular, $0$ and the
one-dimensional faces are always smooth.
\par

{\bf Proposition:}
{\em
If $\tau\leq\sigma$ is a smooth face, then the injections
$\gHom_k(\spann_k K^R_\tau,\,k)\hookrightarrow \bHom(K^R_\tau,k)$ are
even isomorphisms. Moreover,
$\spann_k K^R_\tau =\bigcap_{a^\ai\in\tau}\spann_k K^R_\ai$, and
the latter vector spaces equal
$\spann_k K^R_\ai= M_k$, $(a^\ai)^\bot$, or $0\,$ if $\langle a^\ai,R\rangle
\geq 2$, $=1$, or $\leq 0$, respectively.
}
\par

{\bf Proof:}
Let $f:K_\tau^R\to k$ be quasilinear; it suffices to extend $f$ to a $\Z$-linear map
defined on $\spann_\Z K^R_\tau\subseteq M$.
If $R$ was non-positive on any of the generators of $\tau$, then $K^R_\tau$ would be
empty
anyway. Hence, if (w.l.o.g.) $\tau=\langle a^1,\dots,a^\ak\rangle$, we may assume that
$\langle a^\ai,R\rangle \geq 2$ for $\ai=1,\dots,\al$ and
$\langle a^\aj,R\rangle = 1$ for $\aj=\al+1,\dots,\ak$.\\
$K^R_\tau$ contains the easy part $\tau^\bot\cap\kL_+$; it is no problem at all to
extend $f_{|\tau^\bot\cap\kL_+}$ to a $\Z$-linear function defined on $\tau^\bot\cap M$. 
In general, we have to show that for elements $s^v\in K^R_\tau$ the value
$\sum_v f(s^v)$ only depends on $s:=\sum_v s^v$, not on the summands themselves. 
(Then, $f(s)$ may be defined as this value.)\\
By smoothness of $\tau$ there exist elements $r^1,\dots,r^\al\in K^R_\tau$ 
such that $\langle a^\ai, r^\aj\rangle =\delta_{\ai\aj}$ 
for $1\leq \ai\leq\ak$ and $1\leq\aj\leq\al$. Hence, quasilinearity of $f$ implies
\[
f(s^v)=\sum_{\ai=1}^\al \langle a^\ai, s^v\rangle f(r^\ai) + f(p^v)
\quad \mbox{ with }\;
p^v:= s^v-\mbox{$\sum_\ai$}\langle a^\ai,s^v\rangle r^\ai \in \tau^\bot\cap M\,.
\]
Summing up yields
\[
\sum_v f(s^v)=\sum_{\ai=1}^\al \langle a^\ai, s\rangle f(r^\ai) + \sum_v f(p^v)
=
\sum_{\ai=1}^\al \langle a^\ai, s\rangle f(r^\ai) +
f\big(s-\mbox{$\sum_\ai$}\langle a^\ai,s\rangle r^\ai\big)\,.
\]
Finally, the second claim follows by
$\bigcap_{a^\ai\in\tau}\spann_k K^R_\ai = \bigcap_{\aj=\al+1}^\ak (a^\aj)^\bot
=\spann_k\big(\tau^\bot; r^1,\dots,r^\al\big)= \spann_k K^R_\tau$.
\hfill$\Box$
\par


\neu{E1l-3}
We turn to the remaining part  of the first level and show the vanishing of
$E_1^{p,\geq 2}$ if $Y_\sigma$ is smooth in codimension $p$:
\par

{\bf Theorem:}
{\em
If $\tau\leq\sigma$ is a smooth face, then $\HA^q(K^R_\tau;k)=0$ for $q\geq 2$.
}

{\bf Proof:}
We proceed by induction on $\dim \tau$, i.e.\ we may assume that the vanishing holds
for all proper faces of $\tau$. Let $r(\tau)$ be an arbitrary element of
$\innt(\sigma^\veee\cap\tau^\bot)\cap M$, i.e.\ $\tau=\sigma\cap [r(\tau)]^\bot$.
Then, via $R_g:=R-g\cdot r(\tau)$ with $g\in\Z$,
one obtains an infinite (if $\tau\neq\sigma$) series
of degrees admitting the following two properties:
\kitem
\item[(i)]
$K^{R_g}_\tau=K^R_\tau$ for any $g\in\Z$ (since $R_g=R$ on $\tau$), and
\item[(ii)]
$K^{R_g}_{\tau^\prime}\neq\emptyset$ implies $\tau^\prime\leq\tau$
for any face $\tau^\prime\leq\sigma$ and $g\gg 0$
(since $\langle a^\aj,R_g\rangle\leq 0$ if $a^\aj\notin\tau$).
\kenditem
In particular, in degree $-R_g$ with $g\gg 0$, the first level of our spectral
sequence is shaped as follows:
\kitem
\item
For $p<\dim\tau$ only $\HA^q(K^R_{\tau^\prime};k)$ with $\tau^\prime\leq\tau$ appear as
summands of $E_1^{p,q}$. By the induction hypothesis they even vanish for $q\geq 2$,
\item
for $p=\dim\tau$ it follows that $E_1^{p,q}=\HA^q(K^R_\tau;k)$, and
\item
all vector spaces $E_1^{p,q}$ vanish beyond the $[p=\dim\tau]$--line.
\kenditem
Hence, the differentials $d_r:E_r^{p,q}\to E_r^{p+r,q-r+1}$ are trivial for
$r\geq 1$, $q\geq 2$, and we obtain
\[
T_A^{q+\dim\tau-1}(-R_g)=\HA^q(K^R_\tau;k)
\quad\mbox{ for }\; g\gg 0\,.
\]
Moreover, under this identification, the multiplication
\[
[\cdot x^{r(\tau)}]: T_A^{q+\dim\tau-1}(-R_g) \to  T_A^{q+\dim\tau-1}(-R_{g+1})
\]
is just the identity map.
On the other hand, we may restrict $T^n_\sigma:=T^n_A$ onto the {\em smooth}, open subset
$Y_\tau:=\Spec\, k[\tau^\veee\cap M]\subseteq Y_\sigma$. Since $k[\tau^\veee\cap M]$
equals the localization of $k[\sigma^\veee\cap M]$ by the element $x^{r(\tau)}$, we
obtain
\[
T^n_\sigma \otimes_{k[\sigma^\veee\cap M]}
k[\sigma^\veee\cap M]_{x^{r(\tau)}} =T^n_\tau=0\quad
\mbox{ for }\; n\geq 1.
\]
In particular, any element of $T_A^{q+\dim\tau-1}(-R_g) \subseteq
T_A^{q+\dim\tau-1}$ will be killed by some power of $x^{r(\tau)}$; but this means
$\HA^q(K^R_\tau;k)=0$.
\hfill$\Box$
\par


\neu{E1l-4}
{\bf Corollary:}
{\em
$E_1^{0,q}=E_1^{1,q}=0\,$ for $q\geq 2$.
}
\par

%
%
\section{Applications}\label{app}


\neu{app-1}
The main ingredient for describing the Andr\'{e}-Quillen cohomology
of $Y_\sigma$ will be the complex
$\bHom(K^R_{\kss \bullet},k) = (E_1^{{\kss \bullet},1},\,d_1)$
which is built from
the vector spaces $\bHom(K^R_\tau,k)$ as $C^q(K^R_{\kss \bullet};k)$
was from $C^q(K^R_\tau;k)$ in \zitat{sps}{2}.
It contains $(\spann_k K^R_{\kss \bullet})^\ast$ as a subcomplex.
\par


\neu{app-2}
First, we discuss the case of an {\em isolated singularity} $Y_\sigma$.
Here, our spectral sequence degenerates completely.
\par

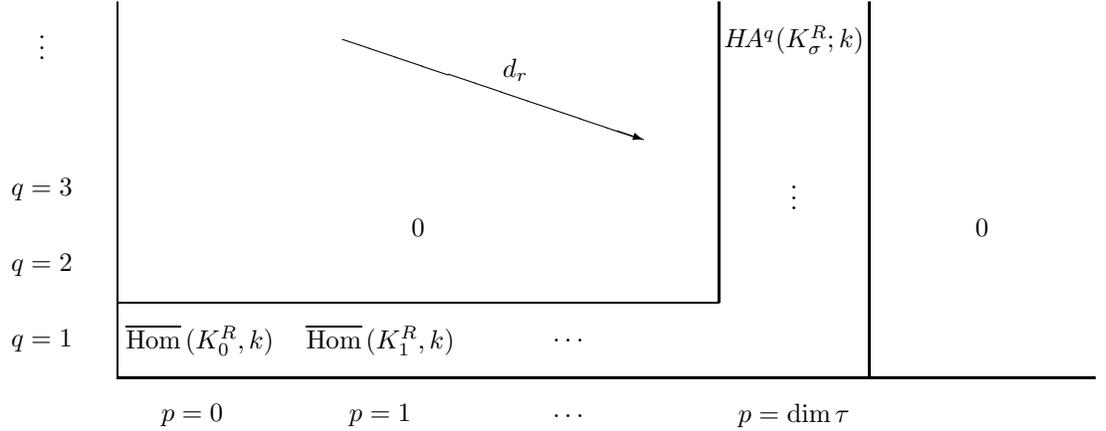
\begin{center}
\unitlength=1.0mm
\linethickness{0.4pt}
\begin{picture}(150.00,60.00)
\put(20.00,10.00){\line(0,1){50.00}}
\put(20.00,10.00){\line(1,0){130.00}}
\put(120.00,20.00){\line(0,1){40.00}}
\put(100.00,60.00){\line(0,-1){40.00}}
\put(120.00,20.00){\line(0,-1){10.00}}
\put(20.00,20.00){\line(1,0){80.00}}
\put(10.00,15.00){\makebox(0,0)[cc]{$q=1$}}
\put(10.00,25.00){\makebox(0,0)[cc]{$q=2$}}
\put(10.00,35.00){\makebox(0,0)[cc]{$q=3$}}
\put(10.00,55.00){\makebox(0,0)[cc]{$\vdots$}}
\put(30.00,5.00){\makebox(0,0)[cc]{$p=0$}}
\put(31.00,15.00){\makebox(0,0)[cc]{$\bHom(K^R_0,k)$}}
\put(55.00,5.00){\makebox(0,0)[cc]{$p=1$}}
\put(55.00,15.00){\makebox(0,0)[cc]{$\bHom(K^R_1,k)$}}
\put(110.00,5.00){\makebox(0,0)[cc]{$p=\dim\tau$}}
\put(110.00,35.00){\makebox(0,0)[cc]{$\vdots$}}
\put(80.00,5.00){\makebox(0,0)[cc]{$\dots$}}
\put(80.00,15.00){\makebox(0,0)[cc]{$\dots$}}
\put(110.00,55.00){\makebox(0,0)[cc]{$\HA^q(K^R_\sigma;k)$}}
\put(60.00,30.00){\makebox(0,0)[cc]{$0$}}
\put(135.00,30.00){\makebox(0,0)[cc]{$0$}}
\put(50.00,55.00){\vector(3,-1){40.00}}
\put(73.00,51.00){\makebox(0,0)[cc]{$d_r$}}
\end{picture}
\vspace{-5ex}
\end{center}
{\bf Proposition:}
{\em
Let $Y_\sigma$ be an isolated singularity. Then, the Andr\'{e}-Quillen
cohomology in degree $-R$ equals
\vspace{-2ex}
\[
T^n_A(-R)=
\renewcommand{\arraystretch}{1.3}
\left\{ \begin{array}{ll}
H^n\big(\bHom(K^R_{\kss \bullet},k)\big)=
H^n\big((\spann_k K^R_{\kss \bullet})^\ast\big) &
\mbox{for }\, 0\leq n\leq \dim\sigma-1\\
H^{\dim \sigma}\big(\bHom(K^R_{\kss \bullet},k)\big) &
\mbox{for }\, n= \dim\sigma\\
\HA^{n-\dim\sigma +1}(K^R_\sigma;k) &
\mbox{for }\, n\geq \dim\sigma+1\,.
\end{array}\right.
\vspace{-2ex}
\]
}
\par

{\bf Proof:}
The first level of the spectral sequence is non-trivial only in
$E_1^{p,1}$ with $0\leq p\leq \dim\sigma$ and
$E_1^{\dim\sigma,q}$ with $q\geq 1$, respectively. Moreover, the complexes
$\bHom(K^R_{\kss \bullet},k)$ and $(\spann_k K^R_{\kss \bullet})^\ast$
are equal up to position $p=\dim\sigma-1$.
\hfill$\Box$
\par


\neu{app-3}
In the {\em general case}, we still have enough information to determine
the deformation relevant modules $T^0_A$, $T^1_A$, and $T^2_A$.
Under the additional hypothesis of smoothness in codimension two,
these results have already been obtained in \cite{T2} with a different proof.
\par

{\bf Proposition:}
{\em
Let $\sigma$ be an arbitrary rational, polyhedral cone with apex in $0$.
Then, for every $R\in M$,
\vspace{-0.5ex}
\[
T^n_A(-R) = H^n\big(\bHom(K^R_{\kss \bullet},k)\big)
\quad\mbox{ for } n=0,1,2.
\]
Moreover, for $n=0,1$, this vector space equals
$H^n\big((\spann_k K^R_{\kss \bullet})^\ast\big)$, too.
}
\par

{\bf Proof:}
This is a direct consequence of Corollary \zitat{E1l}{4}.
\hfill$\Box$
\par


\neu{app-4}
Since the complex $(\spann_k K^R_{\kss \bullet})^\ast$ is much easier
to handle than $\bHom(K^R_{\kss \bullet},k)$, it pays to look for
sufficient conditions for
$T^2_A= H^2\big((\spann_k K^R_{\kss \bullet})^\ast\big)$ to hold.
\par

{\bf Proposition:}
{\em
If $Y_\sigma$ is Gorenstein in codimension two,
i.e.\ for every two-dimensional face $\langle a^\ai,a^\aj\rangle\le\sigma$
there is an element $r(\ai,\aj)\in M$ such that
$\langle a^\ai,r(\ai,\aj)\rangle = \langle a^\aj,r(\ai,\aj)\rangle = 1$,
then
\[
T^2_A= H^2\big((\spann_k K^R_{\kss \bullet})^\ast\big).
\vspace{-3ex}
\]
}
\par

{\bf Proof:}
For edges $\langle a^\ai,a^\aj\rangle\leq\sigma$ we have to show that
$(\spann_k K^R_{\ai\aj})^\ast\hookrightarrow \bHom(K^R_{\ai\aj},k)$ is an
isomorphism. We will adapt the proof of \zitat{E1l}{2}.\\
It may be assumed that $\langle a^\ai,R\rangle; \langle a^\aj,R\rangle \geq 1$
and $r(\ai,\aj)\in K^R_{\ai\aj}$. Let $d:=|\det (a^\ai,a^\aj)|\in\Z$; it is the
smallest positive value of $a^\ai$ possible on elements of
$\kL\cap (a^\aj)^\bot$. We choose an $r^\ai\in\kL\cap (a^\aj)^\bot$
with $\langle a^\ai,r^\ai\rangle=d$; together with $r(\ai,\aj)$ it will play the
same role as $r^1,\dots,r^\al$ did in \zitat{E1l}{2}.
\vspace{0.5ex}\\
{\em Case 1:}\quad
$\langle a^\ai,R\rangle >d$ and
$\langle a^\ai,R\rangle\geq\langle a^\aj,R\rangle$
(in particular, $r^\ai\in K^R_{\ai\aj}$):
Then, elements $s^v$ or $s$ (cf.\ \zitat{E1l}{2}) may be represented as
\[
s=\langle a^\aj,s\rangle\cdot r(\ai,\aj) +
\frac{\kd \langle a^\ai,R\rangle -\langle a^\aj,R\rangle}{\kd d}\, r^\ai
+ \big[ (a^\ai,a^\aj)^\bot-\mbox{elements}\big]\,.
\]
The difference $\langle a^\ai,R\rangle -\langle a^\aj,R\rangle$ is always
divisible by $d$, i.e.\ the coefficients are integers.
\vspace{0.5ex}\\
{\em Case 2:}\quad
$\langle a^\ai,R\rangle; \langle a^\aj,R\rangle\leq d$:
This implies $K^R_{\ai\aj}\subseteq (a^\ai,a^\aj)^\bot +\Z\cdot r(\ai,\aj)
= (a^\ai-a^\aj)^\bot$. In particular, we may use the representation
$s=\langle a^\ai,s\rangle\cdot r(\ai,\aj) +
\big[ (a^\ai,a^\aj)^\bot-\mbox{elements}\big]\,.$
\vspace{1ex}
\hfill$\Box$
\par

{\bf Corollary:}
{\em
Let $Y_\sigma$ be a three-dimensional, toric Gorenstein singularity,
i.e.\ $\sigma=\langle a^1,\dots,a^\an\rangle$ with $a^{\an+1}:=a^1$
being the cone over a lattice polygon embedded in height one.
Then
\[
T^2_A(-R)^\ast=\;^{\kd \mbox{$\bigcap$}_{\ai=1}^\an
\big(\spann_k K^R_{\ai,\ai+1}\big)}
\!\Big/
_{\kd \spann_k \big(\mbox{$\bigcap$}_{\ai=1}^\an K^R_{\ai,\ai+1}\big)}\,.
\]
}
\par


\neu{app-5}
Finally, we would like to mention an alternative to the complexes
$(\spann_k K^R_{\kss \bullet})^\ast$ and $\bHom(K^R_{\kss \bullet},k)$.
Let $E\subseteq\kL_+$ be the (finite) set of non-splittable elements in $\kL_+$;
it is the minimal generator set of the monoid $\kL=\sigma^\veee\cap M$
and gives rise
to a canonical surjection $\pi:\Z^E\to M$.
The relations among $E$-elements are gathered in the $\Z$-module
$L(E):=\ker \pi$.
Every $\kq\in L(E)$ splits into a difference
$\kq=\kq^+ -\kq^-$ with $\kq^+,\kq^-\in \N^E$ and $\sum_v \kq^+_v\,\kq^-_v=0$. We
denote by $\bar{\kq}\in\kL$ the image $\bar{\kq}:=\pi(\kq^+)=\pi(\kq^-)$.
\par

{\bf Definition:}
With $E^R_\tau:=E\cap K^R_\tau$, we define
$L(E^R_\tau):=L(E)\cap \Z^{E^R_\tau}$ and $\bL(E^R_\tau)\subseteq L(E^R_\tau)$
to be the submodule generated by the relations $\kq\in L(E)$ such that $\bar{\kq}\in
K^R_\tau$. (Notice that with $E^R_0=E$ the notation differs slightly from
that in \cite{T2}.)
\par

As usual (cf.\ \zitat{sps}{2}), we may construct complexes from these
finitely generated, free abelian groups. They fit into the following commutative
diagram with exact rows:
\[
\dgARROWLENGTH=0.8em
\begin{diagram}
\node{0}
\arrow{e}
\node{(\spann_k K^R_{\kss \bullet})^\ast}
\arrow{e}
\arrow{s}
\node{k^{E^R_{\kss \bullet}}}
\arrow{e}
\arrow{s,r}{\sim}
\node{\gHom_{\Z}(L(E^R_{\kss \bullet}),k)}
\arrow{e}
\arrow{s}
\node{0}\\
\node{0}
\arrow{e}
\node{\bHom( K^R_{\kss \bullet},k)}
\arrow{e}
\node{k^{E^R_{\kss \bullet}}}
\arrow{e}
\node{\gHom_{\Z}(\bL(E^R_{\kss \bullet}),k)}
\end{diagram}
\]
The cokernel of the second row is $\gExt^1_\Z(L/\bL,k)$. It vanishes
if $\kar\, k=0$.
\par

{\bf Proposition:}
{\em
If $\kar\, k=0$, then the $k$-duals of $T^n_A(-R)$ equal
\[
T^1_A(-R)^\ast=L_k\big(\mbox{$\bigcup$}_\ai E_\ai^R\big)
\big/\mbox{$\sum$}_\ai L_k(E^R_\ai)
\vspace{-2ex}
\]
and
\[
T^2_A(-R)^\ast=\,\ker\big(\oplus_\ai L_k(E_\ai^R) \to L(E)\big)\Big/\,
\im\big(\oplus_{\ai\aj} \bL_k(E_{\ai\aj}^R)\to \oplus_\ai L_k(E_\ai^R)\big)\,.
\vspace{-2ex}
\]
}
\par

{\bf Proof:}
As in the proof of Lemma \zitat{sps}{2}, one obtains that the complex
$k^{E^R_{\kss \bullet}}$ has no cohomology except
$H^0=k^{E\setminus\cup_\ai E^R_\ai}$. Hence, the long exact sequence
for the second row in the above diagram yields isomorphisms
\[
H^{p-1}\big(\gHom_\Z(\bL(E^R_{\kss\bullet}),k)\big)
\stackrel{\sim}{\longrightarrow}
H^{p}\big(\bHom(K^R_{\kss\bullet},k)\big)\quad \mbox{ for }\;p\geq 2.
\]
Taking a closer look at the first terms shows that the
same result is true
for $p=1$ if $E^R_0$ is replaced by $\bigcup_\ai E^R_\ai$.
Finally, we know that $\bL(E^R_\tau)=L(E^R_\tau)$ if $\dim\tau\leq 1$.
\hfill$\Box$
\par


\neu{app-6}
{\bf Remark:}
The vector space $T^1_A(-R)$ has also a convex-geometric interpretation;
it is related to the set of Minkowski summands
of the cross cut $\sigma\cap [R=1]$. For details, we refer to \cite{Flip}.
\par

%
%


{\small
\parbox{9cm}{
Klaus Altmann\\
Institut f\"ur reine Mathematik der\\
Humboldt-Universit\"at zu Berlin\\
Ziegelstr.~13A\\
D-10099 Berlin, Germany\\
E-mail: altmann@mathematik.hu-berlin.de}
\parbox{6cm}{
Arne~B.~Sletsj\o{}e\\
Department of Mathematics\\
University of Oslo\\
Pb.~1053 Blindern\\
N-0316 Oslo, Norway\\
E-mail: arnebs@math.uio.no}}

\end{document}